\documentclass[a4paper, 12pt,  leqno]{article}
\usepackage{amsmath}
\usepackage{amssymb}
\usepackage{latexsym}
\usepackage{eigen}
\newfont{\gothic}{cmr10 at 12pt}

\newcommand{\Section}[2]{\setcounter{equation}{0}
\section[#1]{#2}}
\begin{document}
\title{The First $L^2$-Betti Number of Classifying Spaces for Variations of
Hodge Structures}
\author{
\vspace*{4cm}
\noindent
\parbox[t]{5.5cm}{J. Jost\\
Max Planck Institute \\
for Mathematics in the \\
Sciences \\
Leipzig}
\and
\parbox{1.8cm}{and}
\and
\parbox[t]{3.9cm}{Y. L. Xin\\
Institute of Math.\\
Fudan University \\
Shanghai}}
\date{}
\maketitle

\Section{Introduction}{Introduction}

Classical Hodge theory gives a decomposition of the complex
cohomology of a compact K\"ahler manifold $M$, which carries the
standard Hodge structure $\{H^{p,q}(M),$ $p+q=k\}$ of weight $k$. Deformations of $M$ then lead to variations of the Hodge structure. This is best understood when reformulating the Hodge decomposition  in an abstract manner.
Let $H_C=H_R\oplus C$ be a complex vector space with a real
structure. A Hodge structure on $H_C$ is a decomposition

\beas H_C=\oplus H^{p,q}\quad\text{with}\quad
H^{p,q}=\overline{H^{q,p}}, \quad  p+q=k\eeas Furthermore, 
there exists a $(-1)^k$-symmetric bilinear form $S:H_C\times H_C\to R$
satisfying certain properties. In terms of these abstract structures,  classifying
spaces $D$ for variations of Hodge structures can be
defined.
The subject was initiated by P.~A. Griffiths \cite
{G}. He found that similarly to the Siegel upper-half space from the period matrices of
algebraic curves, certain non-compact homogeneous complex manifolds
arise naturally from the period matrices of general algebraic
varieties.   Those manifolds are also  classifying spaces.

We know that $D=G/V$, where $G$ is a semi-simple Lie group and $V$
is a certain compact  subgroup. $G$ is called a group of Hodge
type. As $V$ typically is a proper subgroup of the maximal compact subgroup $K$ of $G$, $D$ is usually not a symmetric space, but only a homogeneous one. The properties of the projection from $D=G/V$ to $G/K$ were studied by Griffiths and Schmid\cite{GS}. In fact, if one does not attach a marking to the underlying manifold $M$, that is, fix a homology basis, then instead of $G/V$, one rather needs to look at the quotient of that space by the action of a lattice, that is, a discrete subgroup $\Gamma$ of $G$ for which the quotient of $G/V$ by $\Gamma$ has finite volume. Usually, however, that quotient $\Gamma \backslash G/V$ is not compact.

The Hodge decomposition may be interpreted as an action of the
real algebraic group $C^*$ on the cohomology group $H^k(M, C)$
given by
$$z\circ\om=z^p\bar z^q\om$$
for $z\in C^*$ and $\om\in H^{p,q}(M, C)$.

Non-Abelian analogues of the Hodge structure are given by the
action of $C^*$ on moduli spaces of flat bundles. By an important result of Simpson\cite{Sim}, non-Abelian Hodge theory is also related to the
classifying spaces for variations of Hodge structures.

For understanding variations of Hodge structures, it is then important to study the topology of the spaces $G/V$ and $\Gamma \backslash G/V$. The present paper does so in the context of $L^2$-cohomology, by extending previous results for the case of $\Gamma \backslash G/K$, $K$ being a maximal compact subgroup of $G$ as above.

$L^2$-cohomology is the appropriate extension of Hodge theory for harmonic forms to the 
case of a noncompact (complete) manifold $X$ inasmuch as here every $L^2$-cohomology class can be
represented by an $L^2$-harmonic form, see \cite{A}. Besides
offering the possibility to extend Hodge theory to the noncompact
case,  it can be used to obtain topological information about
compact quotients of $X$ by the $L^2$- index theorem   of Atiyah
\cite{A}. This is based on the fact that the operation of a discrete group $\Gamma$ by isometries on $X$ commutes with the Laplacian. Therefore, spaces of harmonic $k$-forms become $\Gamma$-modules, and by constructions from the theory of von-Neumann algebras, they can be assigned dimensions $B^k_\Gamma(X)$. While these dimensions need not be integers in general, the corresponding Euler characteristic is and coincides with the standard one of the quotient $\Gamma \backslash X$ (assuming certain natural assumptions so that the latter is defined). This is the content of Atiyah's theorem. For a detailed discussion, we refer to  the comprehensive volume on
$L^2$-cohomology by W.~L\"uck \cite{L}.

It turns out that $L^2$-cohomology is useful for studying a conjecture of Hopf. Let $\bar M^{2m}$ be a compact manifold of dimension $2m$ with
negative sectional curvature. Hopf conjectured

 \bea (-1)^m \chi(\overline{M}^{2m})
>0. \eea
Dodziuk \cite{D} and Singer \cite{S} suggested to use
$L^2$-cohomology to approach this problem as follows. Take the 
universal covering $M\to \bar M^{2m}$ and show \bea
\mathcal{H}^q(M) = \{0\} \quad \mbox{ for } q \neq m \eea and \bea
\mathcal{H}^m(M) \neq \{0 \}. \eea.
In other words, in contrast to ordinary cohomology, $L^2$-cohomology should concentrate in the middle dimension (for spaces of negative curvature). So far, this approach has given
 partial  answers to the Hopf conjecture, and it is fair to say that these represent the best attack on the problem to date. More precisely, a verification of the conjecture of Dodziuk and Singer has been possible in the following cases:

$\ast$ Symmetric spaces of noncompact type and rank one by
A.~Borel \cite{B}.
\bigskip

$\ast$ K\"ahler hyperbolic manifolds (including quotients of
Hermitian symmetric spaces) by M.~Gromov \cite{G}, extended to the K\"ahler non-elliptic  case (that is, allowing also zero curvature (in which case of course only the first part of the Dodziuk-Singer conjecture can hold) by Jost-Zuo \cite {J-Z} and
Cao-Xavier \cite{C-X}.
\bigskip

$\ast$ For negatively pinched manifolds, Jost-Xin \cite {J-X}
improved the previous results of Donnelly-Xavier \cite{D-X}.
\bigskip

The present paper is devoted to prove

\bete Let $N$ be a classifying space for variations of Hodge structures, $\Gamma$ a lattice on $N$. Then
\beas B_\Ga^1(N)=0\eeas \ete

This result is conceptually different from the ones just quoted because in general $G/V$ is not a space of non-positive sectional curvature since the compact group $K$ and therefore also its (non-trivial) quotient $K/V$ carry some positive curvature. 
The proof will be accomplished by showing that any $L^2$-harmonic 1-form vanishes. An intermediate result will be that the squared norm of any such 1-form is horizontal for the natural Riemannian submersion $G/V \rightarrow G/K$ and only depends on the symmetric space $G/K$. In that sense, the fiber $K/V$ with its positive curvature disappears from the picture, and the situation is reduced to the one of non-positive curvature after all.

\Section{Preliminaries}{Preliminaries} Let $G$ be a semi-simple
Lie group all of whose simple factors are non-compact. Let \go{g}
be its Lie algebra of all left invariant vector fields on $G$ and
$K \subset G$ the Lie subgroup of $G$ whose image in the adjoint
group ad $ G$ is a maximal compact subgroup of ad $G$. Let \go{k}
be the subalgebra of \go{g} corresponding to $K$ and \go{m} the
orthogonal complement of \go{k} in \go{g} with respect to the
Killing form $B(X,Y)$ of \go{g}. Then \bea \mgo{g}= \mgo{m} +
\mgo{k}, \quad [ \mgo{k}, \mgo{k} ] \subset \mgo{k}, \quad
[\mgo{m}, \mgo{m}] \subset \mgo{k}, \quad [\mgo{k}, \mgo{m} ]
\subset \mgo{m}. \eea It is known that the restriction of $B$ to
\go{m} (resp. \go{k}) defines a positive (resp. negative) definite
bilinear form on \go{m} (resp. \go{k}). Hence we can choose a base
$\{X, \hdots, X_n\} $ of \go{m} and a base $\{X_{n+1}, \hdots,
X_{n+r}\} $ of \go{k} with \bea
B(X_i, X_j) &=& \del_{ij}, \\
B(X_\al, X_\beta) &=& - \del_{\al \beta}; \nonumber \eea here and
in the sequel we employ the following range of indices \beas
1  \le & i, j,  \hdots & \le n  \\
n+1 \le & \al, \beta, \ga, \hdots & \le n+r \\
1 \le & a,b,c, \hdots & \le n+r. \eeas Let \beas [X_a, X_b] =
c_{ab}^c X_c. \eeas By (2.1), among the structure constants
$c_{ab}^c$, only $c_{\al \beta}^\ga$, $c_{i j}^\al$, $c_{j
\al}^i$, $ c_{\al j}^i$ can be $\neq 0$.

Let $B(X,Y)$ be the Killing form of \go{g}. It is defined by \bea
B_{ab} =B(X_a, X_b) = \trace( \ad X_a \ad X_b) = c_{ae}^f
c_{bf}^e. \eea Multiplying the Jacobi identity \beas c_{ab}^e
c_{ce}^f + c_{ca}^e c_{be}^f + c_{bc}^e c_{ae}^f = 0 \eeas by
$c_{df}^c$ and summing over the index $f$, we have \beas c_{df}^c
c_{ab}^e c_{ce}^f + c_{df}^e c_{ca}^e c_{be}^f +
c_{df}^c c_{bc}^e c_{ae}^f &\\
=- c_{ab}^e B_{de} + c_{df}^c c_{ca}^e c_{be}^f + c_{df}^c
c_{bc}^e c_{ae}^f & =0. \eeas Denoting \bea c_{ab}^e B_{de} =
c_{dab}, \eea we have \beas
c_{dab} & = & c_{df}^c c_{ca}^e c_{be}^f + c_{df}^c c_{bc}^e c_{ae}^f \\
 & = & c_{de}^f c_{fa}^c c_{bc}^e + c_{df}^c c_{bc}^e c_{ae}^f \\
 & = & c_{bc}^e ( c_{de}^f c_{fa}^c + c_{df}^c c_{ae}^f),
\eeas which is anti-symmetric in $a,d$. Hence $c_{abc}$ is
anti-symmetric in all indices.

(2.2), (2.3) and (2.4) give \bea \sum_{\al, k}  c_{\al i k} c_{\al
j k} = \12 \del_{ij} \eea and \bea \sum_{i, j} c_{i \al j} c_{i
\beta j} + \sum_{\ga, \del} c_{\ga \al \del} c_{\ga \beta \del} =
\del_{\al \beta}. \eea Now let $\{ \om ^a\} $ be invariant forms
dual to $\{X_a\}$. We have the Maurer-Cartan equations \bea d
\om^a = -\12 c_{be}^a \om ^b \wg \om^c. \eea

Let us consider the case when the Cartan subgroup $H\subset G$ is
compact. Let $S^1\subset H$. Then the centralizer of $S^1$ in $G$
is denoted by $V$ and $H\subset V\subset K$. We then have a
complex homogeneous space $G/V$. Denote $\dim H = r_1+1, \; \dim
V=r_1+r_2+1$ and assume that  $r_1+r_2+1<r.$ We also agree on the
range of indices \beas n_1&=&n+(r-r_1-r_2-1),\\
s, t &= & n+1 , \hdots,   n_1,  \\
i_1 &=& 1, \hdots,  n_1, \\
\al_1&= & n_1, \hdots  n+r.   \eeas

We define the Riemannian metric $ds^2$ on $G$ by

$$ ds^2 = \sum( \om^a)^2 . $$

This induces the canonical metric on $G/V$ for which 
$$\Pi:G\to G/V$$
is a Riemannian submersion. We see that $\{X_{i_1}\}=\{X_i,
X_s\}$ are horizontal vector fields and $\{X_{\al_1}\}$ are
vertical vector fields on $G.$ $\{\Pi_*X_i, \Pi_*X_s\}\in T\, G/V$
are orthonormal vector fields, and the dual coframe is
$\{\underline\om^i, \underline\om^s\}$. The Riemannian metric on
$G/V$ is
$$ds^2=\sum (\underline\om^i)^2+\sum (\underline\om^s)^s.$$

Let $\om$ be a $L^2$-harmonic $1$-form on $G/V$
$$\om= u_i\underline\om^i + u_s \underline\om^s.$$

By definition \bea d\om(\Pi_*X_{i_1}, \Pi_*X_{i_2})
&=&\left(\n_{\pi_*X_{i_1}}\om\right)\Pi_*X_{i_2}
-\left(\n_{\pi_*X_{i_2}}\om\right)\Pi_*X_{i_1} \nonumber\\
&=& \Pi_*X_{i_1}(u_{i_2})-\Pi_*X_{i_2}(u_{i_1})
    +\om([\Pi_*X_{i_2}, \Pi_*X_{i_1}]) \nonumber\\
&=& \Pi_*X_{i_1}(u_{i_2})-\Pi_*X_{i_2}(u_{i_1})
     + C_{i_2i_1}^a\om(\Pi_*X_a) \nonumber
\eea

and

\bea \del\om &=& - \left(\n_{\pi_*X_{i_1}}\om\right)\Pi_*X_{i_i}
\nonumber\\
&=&-\Pi_*X_{i_1}(u_{i_1})+\om\left(\n_{\pi_*X_{i_1}}\Pi_*X_{i_1}\right)
\nonumber\\
&=&-\Pi_*X_{i_1}(u_{i_1})+\om(C_{i_1i_1}^{j_1}\Pi_*X_{j_1})
\nonumber\\
&=&-\Pi_*X_{i_1}(u_{i_1}). \nonumber \eea

We then have \bea\Pi_*X_{i_1}(u_{i_2})-\Pi_*X_{i_2}(u_{i_1})
      = C_{i_1i_2}^{j_1}u_{j_1}\eea
and \bea\Pi_*X_{i_1}(u_{i_1})=0.\eea

\Section {Proof of the main theorem}{Proof of the main theorem}

Let $\frak L_{\Pi_*X_s}$ be the Lie derivative with respect to the 
vector field $\Pi_*X_s,$ $g$ the metric tensor on $G/V.$ Then

\beas \left(\frak L_{\Pi_*X_s}g\right)(\Pi_*X_{i_1},\Pi_*X_{j_1})
&=&\Pi_*X_s(\del_{i_1 j_1})\\
&-&g(\frak L_{\Pi_*X_s}\Pi_*X_{i_1}, \Pi_*X_{j_1})
  -g(\Pi_*X_{i_1}, \frak L_{\Pi_*X_s}\Pi_*X_{j_1})\\
&=&-g([\Pi_*X_s, \Pi_*X_{i_1}], \Pi_*X_{j_1})
       -g(\Pi_*X_{i_1}, [\Pi_*X_s, \Pi_*X_{j_1}])\\
&=&-g(\Pi_*C_{si_1}^{k_1}X_{k_1}, \Pi_*X_{j_1})
     -g(\Pi_*X_{i_1}, \Pi_*C_{s j_1}^{k_1}X_{k_1})\\
&=&-C_{s i_1}^{j_1} +C_{s j_1}^{i_1} \eeas which is zero. This
means that $\frak L_{\Pi_* X_s}g=0$ and $\Pi_* X_s$ is a Killing
vector field on $G/V.$ Noting that $\om$ is harmonic and by a result in
\cite{Y-B}
$$\frak L_{\Pi_*X_s}\om=0.$$
By using the formula
$$\frak L_{\Pi_*X_s}\om=(d\circ i_{\Pi_*X_s}+i_{\Pi_*X_s}\circ d)\om,$$
we have
$$d\,u_s = 0,$$
and $u_s$ is constant. $\om$ is an $L^2$-harmonic $1$-form and

\beas | \om |^2 = \sum_i u_i^2 + \sum _s u_s^2, \eeas

is integrable. It turns out that each constant $u_s$ must be zero.
Furthermore, \beas
\Pi_*X_s \sum u_i^2 & = & 2 u_i \Pi_*X_s u_i \\
 & = & 2 u_i (\Pi_*X_i u_s + c_{s i}^j u_j ) \\
 & = & 2 c_{s i }^j u_i u_j = 2 c_{j s i} u_i u_j =0,
\eeas noting that the $c_{abc}$ are anti-symmetric in all indices.

We see that $\Pi_1:G/V\to G/K$ is also a Riemannian submersion
whose fiber is a compact submanifold $K/V$ in $G/V$. Notice that
$K/V$ is a totally geodesic submanifold in $G/V$.

In summary, we have shown that the squared norm of any
$L^2$-harmonic $1$-form $\om$ is horizontal  and  depends only
on the symmetric space $G/K.$

Let $\om$ be a $L^2$-harmonic $1-$form in $N$. We have the
Riemannian submersion \beas \Pi_1: N\to M \eeas with  totally
geodesic fibers, where $M$ is the corresponding symmetric space.
Furthermore, $|\om|^2$ is horizontal. By $L^2$-Hodge theory it is
only necessary to prove that $\om$ vanishes.

Take any unit vector field $n$ in $M$ and any function $b$ in $M$.
We have the horizontal lift $\bar X$ of $bn$. $\bar X$ is the
normal vector field of the fiber submanifold whose length is
constant along the fibers. Choose an orthonormal frame field $\{
e_i \}$ in $ M$, and call its horizontal lift $\{ \bar {e_i} \}$.
$\{\bar{ e_s} \}$ is an orthonormal frame on the fiber. Thus, $\{
\bar{e_i}, \bar{e_s}\}$ is an orthonormal frame field on $N$.
Therefore $\lag \n_{\bar{e_s}} \bar X, \bar{e_s} \rag $ is a
multiple of the mean curvature with respect to the normal
direction $n$. It is zero since the fibers are totally geodesic.
$\Div \bar X$ can be computed in the base manifold $M$. We also
have

\beas \lag \om \odot \om , \n \bar X \rag = u_i u_j \lag \n_{e_i}
X, e_j \rag. \eeas  Hence \beas \lag S_\om, \n X \rag = \12 | \om
|^2 \Div X - \lag \om \odot \om, \n X \rag \eeas can be computed
in $M$, provided $X$ is of the above type, where $S_\om$ is the
stress-energy tensor of $\om$.

Choose $D= B_R(x_0)$, a geodesic ball in $M$  of 
radius  $R$  with  center  in $x_0 \in M$.
Its  boundary  is a  geodesic
sphere  $S_R(x_0)$ in  $M$. Let \\
$\bar D= \Pi_1^{-1}(D) \subset N$. $\prt\bar D$ is compact since
the fiber is compact. Let $X = r \dr$, which is a smooth vector
field in $M$. Let $\bar X$ be the horizontal lift of $X$. Since
the fiber submanifold is orthogonal to the horizontal vector
field, $\bar X$ is also a normal vector field on $\prt\bar D$. Its
length is equal to $r$. Thus, for any $L^2$-harmonic $1$-form
$\om$ in $N$, we  have

\bea \int_{\prt\bar D} \12 |\om|^2 \lag\bar X, n \rag \ast 1 -
\int_{\prt\bar D} \lag i_{\bar X} \om, i_n \om \rag \ast 1 \\
= \int_{\prt\bar D} \12 R|\om|^2 \ast 1 - \int_{\prt\bar D} R \lag
i_{\dr} \om, i_{\dr} \om \rag \le \12 R \int_{\prt\bar D} | \om
|^2 \ast 1. \nonumber \eea

On the other hand,

\beas \n_\frac{\prt}{\prt r} X &=& \frac{\prt}{\prt r}, \;
\n_{e_{s'}} X = r \He(r) (e_{s'}, e_{t'}) e_{t'}, \\
\Div X &=&1 + r \He(r)(e_{s'}, e_{s'}), \eeas where $\{e_i \} = \{
e_{s'}, \frac{\prt}{\prt r} \}\;(s',t'=1, \hdots, n-1)$ is an
orthonormal frame field in $D$.  Therefore \beas \lag \om \odot
\om, \n X \rag = \big| i_\frac{\prt}{\prt r} \om\big|^2 +\lag
i_{e_{s'}} \om, i_{e_{t'}} \om \rag r \He(r)(e_{s'}, e_{t'}),
\eeas and hence

\beas \lag S_{\om}, \n X \rag &=& \12 |\om|^2 (1 + r \,
\mbox{Hess}\, (r) (e_{s'}, e_{s'}))
-\big| i_\frac{\prt}{\prt r} \om \big|^2 \\
&-& \lag i_{e_{s'}} \om, i_{e_{t'}} \om \rag r \He(r)(e_{s'},
e_{t'}) \nonumber\\& = & \left( \12 \sum_{s'} r \He (r) ( e_{s'},
e_{s'})-\12\right) \left| i_{\dr} \om \right|^2 \\& & +
\sum_{s'}\left( \12 + \12 \sum_{t'} r \He (r) (e_{t'},
e_{t'})\right) \lag i_{e_{s'}} \om, i_{e_{s'}} \om
\rag \\
& & - \sum_{s',t'} r \He (r) (e_{s'}, e_{t'}) \lag i_{e_{s'}} \om,
i_{e_{t'}} \om \rag. \eeas

Choose a local orthonormal frame field $\{ e_{s'}\}$ near $x$ in
$S_r(x_0)$, such that $\He(r)$ is diagonalized at $x$. By parallel
translating along the radial geodesics from $x_0$ we have a local
orthonormal frame field in $M$. We have at $x$

\bea &&\lag S_\om , \n X \rag  = \left(\12 \sum_{s'} r \He (r)
(e_{s'}, e_{s'}) - \12 \right) \left| i_\frac{\prt}{\prt r} \om
\right|^2  \\
&+& \sum_{s'} \left( \12 + \12 \sum_{t'} r \He (r) (e_{t'},
e_{t'}) - r \He (r) (e_{s'}, e_{s'}) \right) \lag i_{e_{s'}} \om,
i_{e_{s'}} \om \rag.\nonumber\eea

First of all, by the Hessian comparison theorem \bea \12 \sum_s r
\He (r) (e_s, e_s) - \12 \ge \frac{n-2}{2} >0. \eea To estimate
the coefficients of the second term of (3.2) let \beas A_s =
\sum_{t'} \He (r) (e_{t'}, e_{t'}) - 2 \He (r) (e_{s'}, e_{s'}).
\eeas Since \beas \n_\frac{\prt}{\prt r} \He (r) (e_{s'}, e_{s'})
& = &
\lag \n_\frac{\prt}{\prt r} \n_{e_{s'}} \frac{\prt}{\prt r}, e_{s'} \rag \\
& = & - \lag R(\frac{\prt}{\prt r}, e_{s'}) \frac{\prt}{\prt r},
e_{s'} \rag +
\lag \n_{ [\frac{\prt}{\prt r}, e_{s'}]}  \frac{\prt}{\prt r}, e_{s'} \rag \\
& = & - \lag R ( \frac{\prt}{\prt r}, e_{s'}) \frac{\prt}{\prt r},
e_{s'} \rag - \lag \n_{e_{s'}} \frac{\prt}{\prt r}, \n_{e_{s'}}
\frac{\prt}{\prt r} \rag, \eeas we have \beas \frac{d}{d r}(\Delta
r) = - \Ric ( \frac{\prt}{\prt r}, \frac{\prt}{\prt r} ) - | \He
(r) |^2, \eeas

Moreover, if the sectional curvature of $M$ satisfies $-a^2\le
K\le 0$ and its Ricci curvature Ric $\le -b^2$, then

\beas \frac{d A_s(r)}{d r} & = & (-\Ric(\dr, \dr ) -|\He(r)|^2)\\
&\qquad&+  2 \lag R(\dr, e_{s'} )\dr, e_{s'} \rag
+ 2 \lag \n_{e_{s'}} \dr, \n_{e_{s'}} \dr \rag \\
& \ge & b^2 -2 a^2 +2 \lag \n_{e_s} \dr, \nabla_{e_s} \dr \rag -|\He(r) |^2 \\
& = & b^2 - 2 a^2 + \left( \He (r) (e_{s'}, e_{s'}) \right)^2 \\
&\qquad& - \sum_{u' ,v' \neq s'} \He (r) (e_{u'}, e_{v'}) \He (r)
(e_{u'}, e_{v'}). \eeas

For the classifying spaces $N$, the base manifold $M$ of the
Riemannian submersion \hfill $\Pi_1: N\to M$ \hfill is \hfill
$SO(p, 2q)/SO(p)\times SO(2q)\; (q\ge 2)$ \hfill or \hfill
$Sp(m+n)/Sp(m)\times Sp(n)$. We know that \cite {W}

\begin{table}[h]
\renewcommand{\arraystretch}{1.2}
\caption{}
\begin{tabular*}{13.7cm}[t]{|@{\extracolsep{\fill}}c|c|c|}
\hline
{\bf Type}  & {\bf Sec. Curvature} &{\bf Ricci Curvature} \\
\hline \,$SO(p, 2q)/SO(p)\times SO(2q),\; p>1$  & $ -2 \le K \le
0$ & $-(p+2q-2)$\\
\hline \,$SO(p, 2q)/SO(p)\times SO(2q),\; p=1$  & $  K = -1$ &
$-(p+2q-2)$
\\ \hline\, $Sp(m,n)/Sp(m)\times Sp(n)$ &  $-4 \le K \le 0$ &
$-4(n+m+1)$ \\
$m,n>1$ && \\
\hline $Sp(m,n)/Sp(m)\times Sp(n)$  & $K = -4 $ &
$-4(n+m+1)$ \\
$m=n=1$&&\\
\hline $Sp(m,n)/Sp(m)\times Sp(n),$ & $-4
\le K \le -1$ & $-4(n+m+1)$\\
$\min(m, n)=1$ and $m\ne n$& &\\ \hline
\end{tabular*}
\end{table}
In any case we have

\beas \frac{d A_s(r)}{d r}\ge \left( \He (r) (e_{s'}, e_{s'})
\right)^2 - \sum_{u' ,v' \neq s'} \He (r) (e_{u'}, e_{v'}) \He (r)
(e_{u'}, e_{v'}). \eeas
 Noting  that  the  sectional  curvature of
\hfill $M$ \hfill is \hfill nonpositive \hfill and \hfill each \\
$\He (r) (e_t, e_t) \ge \frac{1}{r} >0$, \bea \frac{d A_s(r) }{d
r} & \ge & \left( \He (r) (e_s, e_s) \right)^2 -
\left( \sum_{t \neq s} \He(r) (e_t, e_t) \right)^2 \nonumber \\
& = & \left( \He (r) (e_s, e_s) - \sum_{t \neq s} \He (r) (e_t,
e_t) \right)
\Del r \\
& = & -A_s(r) \Del r. \nonumber \eea Since $A_s(0) >0$ because the
dimension is at least $4$, we may deduce from (3.4) that $A_s(r)
>0$ for all $r>0$. Altogether, we conclude that \beas \lag S_\om,
\nabla X \rag \ge \mbox{ const. } |\om|^2 \eeas for a positive
constant. If $|\om|\neq 0$ there exists $R_0>0$ such that when
$R\ge R_0$

\bea\int_{\bar D}\left <S_{\om}, \n\bar X\right>*1\ge C>0.\eea

We have the basic inequality for the stress energy tensor \cite{X}

\bea \int_{\prt \bar D} \12 |\om|^2 \lag\bar X, n \rag \ast 1 =
\int_{\bar D} \lag S_\om, \n\bar X \rag \ast 1 + \int_{\prt \bar
D} \lag i_{\bar X} \om, i_n \om \rag \ast 1.\eea

By (3.1),(3.5)and (3.6) we obtain

\beas\int_{\prt \bar D}|\om|^2*1\ge \frac{2C}R\eeas

and

\beas\int_N|\om|^2*1\ge\int_{R_0}^\infty  d\,R\int_{\prt \bar
D}|\om|^2*1=\infty\eeas

which contradicts the $L^2$-assumption on $\om$.

\bigskip

{\noindent\bf Acknowledgements:}  The second author thanks the Max
Planck Institute for Mathematics in the Sciences in Leipzig for
providing good working conditions during the preparation of this
paper, and also 973 project and SFECC for support.

\end{document}